\documentclass[10pt]{article}

\def\pref#1{(\ref{#1})}
\input{psfig.sty}
 \input epsf 
\def\morespace{\vspace{5pt}}
\def\pochhammer#1#2{(#1)_{#2}}
\def\Frac#1#2{\frac
{
 {\raise.6ex
 \hbox{$\displaystyle#1$}}
}
{
 {\lower.6ex
 \hbox{$\displaystyle#2$}}
 }
}
\def\wt{\widetilde}
\def\intp{\int_0^\infty}

\def\bo{{\cal O}}
\def\D{{\cal D}}
\def\P{{\cal P}}
\def\C{{\cal C}}
\def\Vp{{\cal V}}
\def\erfc{{{\rm erfc}}}

\def\ph{{{\rm ph\,}}}
\def\la{\lambda}
\def\tfrac#1#2{{{\lower.6ex
\hbox{$\scriptstyle#1$}}\over 
{\raise.7ex
\hbox{$\scriptstyle#2$}}}}
\def\F#1#2#3#4{{}_2F_1\left(\matrix{#1,#2\cr#3\cr};\,#4\right)}
\def\RR{{{\rm I}\!{\rm R}}}

\begin{document}

 \title{
Parabolic Cylinder Functions:
Examples of Error Bounds For Asymptotic Expansions}

\author{Raimundas Vidunas\\ 
  CWI, P.O. Box 94079, 1090 GB Amsterdam, The Netherlands\\
and \\
  Korteweg-de Vries Instituut voor Wiskunde, University of Amsterdam, \\
Plantage Muidergracht 24, 
1018 TV Amsterdam, 
The Netherlands\\
e-mail: {\tt vidunas@cwi.nl, vidunas@wins.uva.nl}\\
\and
Nico M. Temme\\
  CWI, P.O. Box 94079, 1090 GB Amsterdam, The Netherlands\\
e-mail: {\tt  nicot@cwi.nl}
}
\date{\today}
\maketitle
\begin{abstract}
\noindent
Several asymptotic expansions of parabolic cylinder functions are discussed
and error bounds for remainders in the expansions are presented. In
particular Poincar{\'e}-type expansions for large values of the argument $z$
and uniform expansions for large values of the parameter are considered.
The expansions are based on those given
in \cite{Olver:1959:UAE} and on modifications of these expansions  given in
\cite{Temme:2000:NAP}. Computer algebra techniques are used for obtaining
representations of the bounds and for numerical computations.
 \end{abstract}

\vskip 0.8cm \noindent
1991 Mathematics Subject Classification:
41A60, 33C15, 33C10, 30E15, 33F05, 65D20.
\par\noindent
Keywords \& Phrases:
Parabolic cylinder functions,
asymptotic expansions, 
error bounds for remainders,
Airy functions,
numerical evaluation of special functions.
\section{Introduction}

The solutions of the differential equation
\begin{equation}
\frac{d^2y}{dz^2}-\left(\tfrac14z^2+a\right)y=0 \label{eq:IN.I1}
\end{equation}
are associated with the parabolic cylinder in harmonic analysis; see
\cite{Weber:1869:UIP}. The solutions are called parabolic cylinder functions
and are entire functions  of
$z$. Many properties are given in \cite{Miller:1955:TWP} and 
\cite{Abramowitz:1964:HMF}; for applications to physics and many more
properties see \cite{Buchholz:1969:CHF}.

As in \cite{Miller:1955:TWP} and 
\cite{Abramowitz:1964:HMF}, Chapter 19,
we denote  two standard solutions of \pref{eq:IN.I1} 
by $U(a,z), V(a,z)$. 
Another well-known notation for the parabolic cylinder function is $D_\nu(z)$.
The relation  between $D_\nu(z)$ and $U(a,z)$ is
\begin{equation}\label{eq:IN.I13}
D_\nu(z)=U(-\nu-\tfrac12,z).
\end{equation}
Wronskian relations for the solutions $U(a,z), U(a,-z), V(a,z)$ of 
\pref{eq:IN.I1} are:
\begin{equation}\label{eq:IN.I14}
U(a,z)V'(a,z)-U'(a,z)V(a,z)=\sqrt{{2/\pi}},
\end{equation}
\begin{equation}\label{eq:IN.I15}
U(a,z)\frac{dU(a,-z)}{dz}-U'(a,z)U(a,-z)=\frac{\sqrt{2\pi}}{\Gamma(a+\tfrac12)},
\end{equation}
which shows that $U(a,z)$ and $V(a,z)$ are independent solutions 
of \pref{eq:IN.I1} for all values of $a$. 

Other relations and connection formulae are
\begin{equation}\label{eq:IN.I16}
U(a,z)=\frac{\pi}{\cos^2\pi a\,\Gamma(a+\tfrac12)}
\left[V(a,-z)-\sin\pi a\,V(a,z)\right],
\end{equation}
\begin{equation}\label{eq:IN.I17}
\pi V(a,z)=\Gamma(\tfrac12+a)\left[\sin\pi a\,
U(a,z)+U(a,-z)\right],
\end{equation}
\begin{equation}\label{eq:IN.I18}
\sqrt{2\pi}\,U(-a,iz)
=\Gamma(\tfrac12+a)\left[e^{-i\pi(\frac12a-\frac14)}\,U(a,z)+
e^{i\pi(\frac12a-\frac14)}\,U(a,-z)\right],
\end{equation}
\begin{equation}\label{eq:IN.I19}
U(a,z)
=ie^{\pi ia}U(a,-z)+\frac{\sqrt{2\pi}}{\Gamma(a+\tfrac12)}
e^{\frac12\pi i(a-\frac12)}U(-a,iz),
\end{equation}
\begin{equation}\label{eq:IN.I20}
U(a,z)
=-ie^{-\pi ia}U(a,-z)+\frac{\sqrt{2\pi}}{\Gamma(a+\tfrac12)}
e^{-\frac12\pi i(a-\frac12)}U(-a,-iz).
\end{equation}

In \cite{Olver:1959:UAE}
an extensive collection of asymptotic expansions for the parabolic cylinder
functions as $|a|\to\infty$ has been derived from the differential equation
\pref{eq:IN.I1}. The expansions are valid for complex values of the parameters and are
given in terms of elementary functions and Airy functions.
In \cite{Temme:2000:NAP} modified expansions are given, which have as extra
feature that the expansions are also valid when $a$ is fixed and $z$ is large.
The coefficients of the modified expansions are different from
those of Olver's expansions, and they can be generated by recursion formulas.

When Olver published his results, his later work on bounds for remainders in
asymptotic expansions was not available, and, as he remarked in
\cite{Olver:1975:UPA}, the construction of error bounds for asymptotic
expansions of the parabolic cylinder functions was an important problem to be
considered.
In this paper we discus error bounds for the remainders of
the standard Poincar\'e-type expansions of $U(a,z)$, 
and of some of the uniform expansions.
\section{Poincar\'{e}-Type Expansions}\label{sec:PT}
These expansions are  for large $z$ and $a$ fixed. 
They are given in \cite{Abramowitz:1964:HMF} and 
derived in \cite{Whittaker:1952:CMA}. We have

\begin{equation}\label{eq:PT.E1}
U(a,z)\sim e^{-\frac14z^2}\,z^{-a-\frac12}\,
\sum_{s=0}^\infty
(-1)^s\frac{\pochhammer{a+\tfrac12}{2s}}{s!(2z^2)^s},\quad
|\ph{z}|<\tfrac34\pi,
\end{equation}
\begin{equation}\label{eq:PT.E2}
V(a,z)\sim \sqrt{\tfrac{2}{\pi}}e^{\frac14z^2}z^{a-\frac12}
\sum_{s=0}^\infty
\frac{\pochhammer{-a+\tfrac12}{2s}}{s!(2z^2)^s}, \quad
|\ph{z}|<\tfrac14\pi.
\end{equation}
By using \pref{eq:IN.I19} and  \pref{eq:IN.I20} the sector of validity 
of \pref{eq:PT.E3} can be modified, and compound expansions are 
obtained:
\begin{equation}
\begin{array}{ll}\label{eq:PT.E3}
&{\displaystyle{U(a,z)\sim e^{-\frac14z^2}z^{-a-\frac12}\,
\sum_{s=0}^\infty
(-1)^s\frac{\pochhammer{a+\tfrac12}{2s}}{s!(2z^2)^s}}}\\
&\quad\quad+
{\displaystyle{i\frac{\sqrt{2\pi}}{\Gamma(a+\tfrac12)}e^{-a\pi i}
e^{\frac14z^2}z^{a-\frac12}\,
\sum_{s=0}^\infty
\frac{\pochhammer{-a+\tfrac12}{2s}}{s!(2z^2)^s}}},\quad
\tfrac14\pi<\ph{z}<\tfrac54\pi, 
\end{array}
\end{equation}
\begin{equation}
\begin{array}{ll}\label{eq:PT.E4}
&{\displaystyle{U(a,z)\sim e^{-\frac14z^2}z^{-a-\frac12}\,
\sum_{s=0}^\infty
(-1)^s\frac{\pochhammer{a+\tfrac12}{2s}}{s!(2z^2)^s}}}\\
&\quad\quad
{\displaystyle{-i\frac{\sqrt{2\pi}}{\Gamma(a+\tfrac12)}e^{a\pi i}
e^{\frac14z^2}z^{a-\frac12}\,
\sum_{s=0}^\infty
\frac{\pochhammer{-a+\tfrac12}{2s}}{s!(2z^2)^s}}},\quad
-\tfrac54\pi<\ph{z}<-\tfrac14\pi.
\end{array}
\end{equation}

With these results we can also obtain compound expansions for $V(a,z)$ 
for other sectors than given in \pref{eq:PT.E2}.

\subsection{Error bounds of the expansions.}\label{sec:PT.ER}
Bounds for remainders in the Poincar{\'e}-type expansion follow 
from \cite{Olver:1965:ASW},
where results are given for Whittaker functions. 
The function $U(a,z)$ is a special case of this function.
The relation is
\begin{equation}\label{eq:PT.E7}
U(a,w)= 2^{-\frac12a} w^{-\frac12}
W_{k,m}(z),
\quad k=-\tfrac12a,\quad m=\tfrac14, 
\quad z=\tfrac12w^2.
\end{equation}
The asymptotic expansion for the Whittaker function reads
\begin{equation}\label{eq:PT.E8}
W_{k,m}(z)= z^k e^{-\frac12z}\sum_{s=0}^{n-1} \frac{a_s}{z^s} +\epsilon_n(z), 
\quad a_s=(-1)^s\frac{(a+\frac12)_{2s}}{2^{2s}\,s!},
\quad n=0,1,2,\ldots\ .
\end{equation}
We introduce the following quantities. Let 
\begin{equation}
    \begin{array}{ll}\label{eq:PT.E9}
    &\kappa=|a|, \quad \sigma=\Frac{\kappa}{|z|}, \quad
    \alpha=\Frac1{1-\sigma},\quad
    \beta
    =\tfrac12+\tfrac12\sigma+\tfrac12\sigma(1-\sigma)^{-1}|z|^{-1},\\
    &\delta=|\tfrac14a^2+\tfrac3{16}|+
    \sigma(1+\tfrac14\sigma)(1-\sigma)^{-2},\\
    \end{array}
\end{equation}
assuming that $\sigma<1$.
Then remainder $\epsilon_n(z)$
and its derivative can be bounded as follows
\begin{equation}
	\left.\begin{array}{ll}\label{eq:PT.E10}
        |\epsilon_n(z)|&\le\\
	&\\
	|\beta^{-1}\epsilon_n'(z)|&\le\\
	\end{array}\right\}
	2\alpha\left|z^ke^{-\frac12z}a_n\right| \Vp_\P\left(t^{-n}\right)
	\exp\left[2\alpha\delta  \Vp_\P\left(t^{-1}\right)\right],
\end{equation}

$\Vp(f)$ denotes the variational operator (see
\cite{Olver:1974:ASF}), which for 
continuously differentiable functions in a real interval $[a,b]$
is defined by 
\begin{equation}\label{eq:PT.E11}
\Vp_{a,b}(f)=\int_a^b|f'(x)|\,dx.
\end{equation}
For a holomorphic function $f(z)$ in a complex domain the
variational operator along a smooth arc $\C$ parameterized by
$z(\tau)$, $\alpha<\tau<\beta$ in which $\tau$ is the arc parameter
and $z'(\tau)$ is continuous and nonvanishing in the
closure of $(\alpha,\beta)$, we have
\begin{equation}\label{eq:PT.E12}
\Vp_{\C}(f)=\int_\alpha^\beta\left|f'[z(\tau)]z'(\tau)\right|\,d\tau.
\end{equation}
Along a path $\P$ that is a finite chain of smooth arcs
(of straight lines, for example)  $\Vp_\P$
can be defined as the sum of the contributions from the arcs.

In the bounds given in \pref{eq:PT.E10} the path $\P$ 
links the point $z$ (with $\ph{z}\in(-\frac32\pi,\frac32\pi)$)
to $+\infty$, such that on $\P$ the condition is fulfilled that
$\Re(t+a\ln t)$ is monotonic. 

We consider the bounds in \pref{eq:PT.E10} for $z$
in the sector $[-\pi,\pi]$, that is, for $w$ used in 
\pref{eq:PT.E7} with $\Re w\ge 0$. For other values of $w$ the
relations in \pref{eq:IN.I19} and \pref{eq:IN.I20} 
can be used for computing the function $U(a,z)$.

In \cite{Olver:1965:ASW} simple bounds are given for the variation
$\Vp_\P(t^{-n})$ appearing in \pref{eq:PT.E10}, for $z$ in
certain regions in the complex plane. In Figure \ref{fig:fig1} we show these
regions in the $z-$plane and corresponding regions in the
$w-$plane; we only show the regions in $\Im z\ge0, \Im w\ge 0$,
but they should be extended by including the conjugated parts.

\begin{center}
\epsfxsize=12cm \epsfbox{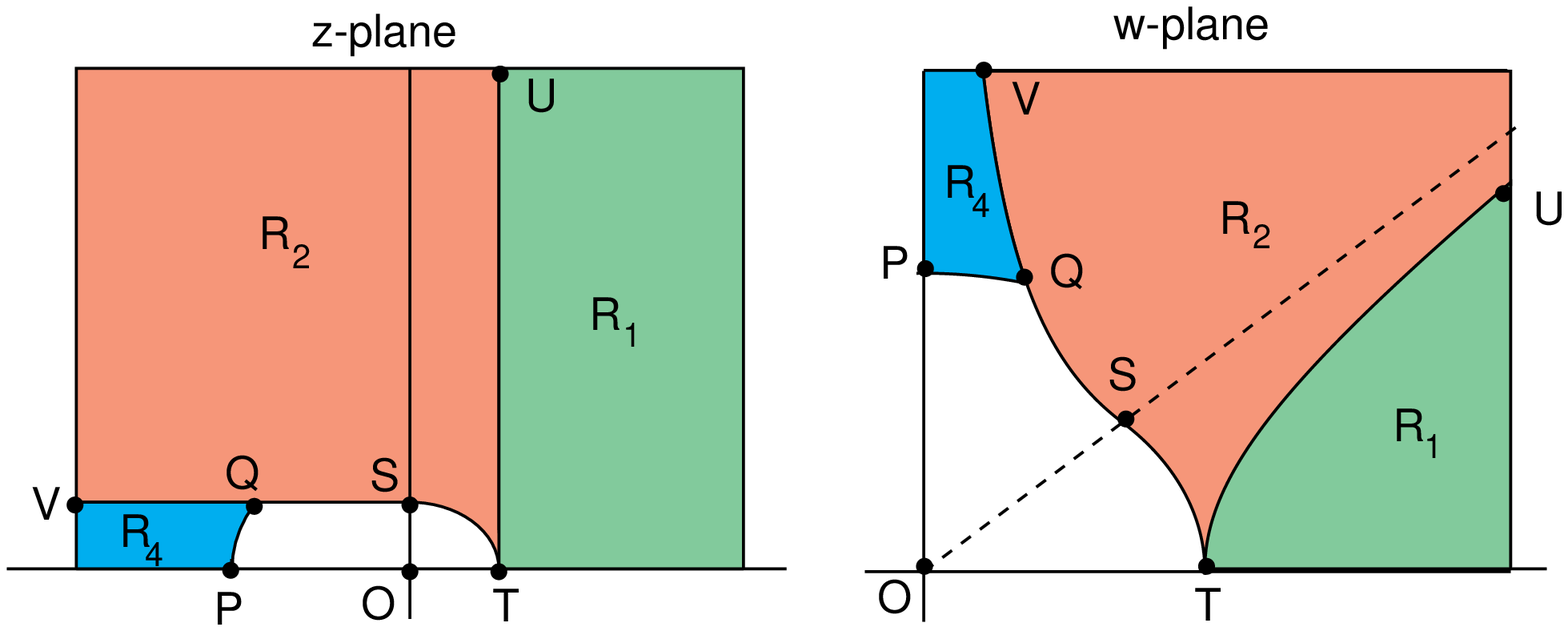}\label{fig:fig1}
\vspace*{0.5cm}

{\bf Figure \ref{fig:fig1}.}\ The regions $R_1, R_2$ and $R_4$ used for 
bounding the variation $\Vp_\P\left(t^{-n}\right)$
appearing in \pref{eq:PT.E10}.
\end{center}
The arc $PQ$ is a circular arc, with radius $2\kappa$; $Q$ is the
point $-\sqrt{3}\kappa+i\kappa$;
$S$ is the point $i\kappa$; the arc $ST$ is a circular arc, 
with radius $\kappa$. As in \pref{eq:PT.E9}
$\kappa=|a|$. 

The region $R_1$ is the half-plane $\Re z\ge \kappa$;
$R_2$ is the region above the curves $VQ, QS, ST$ 
and the conjugated  $z-$values;
$R_4$ is the region with $|z|\ge2\kappa, |\Im z|\le \kappa$. The
corresponding regions in the $w-$plane follow from $z=\frac12w^2$.

The following upper bounds for $\Vp_\P\left(t^{-n}\right)$ are derived 
in \cite{Olver:1965:ASW}:
\begin{equation}
	\begin{array}{ll}\label{eq:PT.E13}
	|z|^{-n}, & z\in R_1,\\
	\chi(n)|z|^{-n}, & z\in R_2,\\
	\left[\chi(n)+\sigma v^2n\right] v^n |z|^{-n}, & z\in R_4,\\
	\end{array}
\end{equation}
where
\begin{equation}\label{eq:PT.E14}
	\chi(n)=\frac{\sqrt{\pi}\Gamma\left(\tfrac12n+1\right)}
	{\Gamma\left(\tfrac12n+\tfrac12\right)},\quad
	v= \left(\tfrac12+\tfrac12\sqrt{1-4\sigma^2}\right)^{-1/2}.
\end{equation}
For $\Vp_\P\left(t^{-1}\right)$ set $n=1$ in \pref{eq:PT.E13}.

For $z\in R_4$, the quantities $\alpha$, $\beta$ and $\delta$ of
\pref{eq:PT.E9} should be modified: replace $\sigma$ by $v\sigma$
and $|z|^{-1}$ by $v|z|^{-1}$.

For the parabolic cylinder function $U(a,w)$ these bounds are applicable
if $|w|$ is large compared with $\sqrt{|a|}$. For example, 
if $w\in R_4$ we need $|w|\ge 2\sqrt{|a|}$. When using the
asymptotic expansion for computations, the restrictions on the bounds
are not unrealistic because the asymptotic expansion does not make sense 
if $w$ is not large compared with $\sqrt{|a|}$.

\subsection{Are these domains optimal?}\label{opt}
We have used the same values as in \cite{Olver:1965:ASW} for showing the
regions $R_1$, $R_2$ and $R_4$. When verifying Olver's analysis we have 
found that the boundaries of the regions $R_1$ and $R_2$ can be modified. 

\subsubsection{Extending region {\protect\boldmath $R_1$}}
For $R_1$ we can use the condition $\Re\, z > \max[0,\Re(-a)]$. 
To verify this condition, let
\begin{equation}\label{eq:PT:v1}
a=u+iv=|a|e^{i\alpha},\quad z=x+iy=re^{i\phi}.
\end{equation}
In Olver's analysis the path in the $z-$plane has to be selected 
along which 
\begin{equation}\label{eq:PT:v2}
F(x,y)=\Re(z+a\ln z)=x+\tfrac12u\ln(x^2+y^2)-v\arctan(y/x)
\end{equation}
is monotonic. Olver chooses an optimal path with the constant 
argument $\phi$. Substituting $y=x\tan\phi$ in $F(x,y)$, we find
$dF/dx=(u+x)/x$. This is positive on the path when  
$\Re\, z > \max[0,\Re(-a)]$. 
One can also use a vertical path $x=x_0, y\ge y_0$
upwards,  if $\phi\in(\alpha,\alpha+\pi)$.  
This follows from 
\begin{equation}\label{eq:PT.v3}
\frac{dF(x_0,y)}{dy}=\frac{uy-vx_0}{x_0^2+y^2}.
\end{equation}
Similarly, one can use the vertical path
downwards if $\phi\in(\alpha-\pi,\alpha)$.
In Olver's approach vertical paths are not used. 

\subsubsection{Extending region {\protect\boldmath $R_2$}}
For $R_2$ the condition $\Im\, z > \max[0,\Im(-a)]$ can be used.
We verify this by taking as Region 2 (as in Olver's approach)
the domain with points $z_0=x_0+iy_0$ from which we can draw a half line 
with  the equation  $x_0 x + y_0 y = x_0^2 + y_0^2$.
This line is perpendicular (at $z=z_0$) to the line
from the origin to $z_0$.  Another equation for the line is
$y = y_0 - \frac{x_0}{y_0} (x-x_0)$.
On this path we have
\begin{equation}\label{eq:PT.v4}
\frac{dF}{dx} = 1 + \frac{u (x_0^2+y_0^2)}{y_0^2 (x^2+y^2)} (x-x_0)
        +\frac{v (x_0^2+y_0^2)}{y_0 (x^2+y^2)},
\end{equation}
and we see that $dF/dx>0$ at $z=z_0$ if $0<-v<y_0$ or $v>0, y_0>0$.
With this conditions $dF/dx>0$ for all $x>x_0$.
This explains that we can extend Region 2 to the domain where 
$\Im z>\max(0,\Im(-a))$.

\subsection{Application to the error function}
We applied these bounds for the case $a=\tfrac12$, 
which corresponds to the error function:
\begin{equation}\label{eq:PT.E15}
W_{-\frac14,\frac14}(z)= \sqrt{\pi}\, z^{\frac14}\,e^{z}\,\erfc\,\sqrt{z}.
\end{equation}
We computed 
\begin{equation}\label{eq:PT.E16}
\rho=\frac{|\epsilon_n(z)|}{\epsilon_n^{(e)}(z)},
\end{equation}
where $\epsilon_n(z)$ is the exact error (see \pref{eq:PT.E8}), 
$\epsilon_n^{(e)}(z)$ 
the estimated error (see \pref{eq:PT.E10}), for several values of $n$ 
and $\theta=\ph{z}$ as given in Table \ref{tab:PT.E17}. We observe that the ratio $\rho$ 
is almost always less than $\frac13$, 
and that in $R_2$, where $\theta=\frac{j}8\pi$
with $j=4, 5, 6, 7$, the estimated error is quite large compared with the
real error.

\begin{table*}
\caption{Ratios $\rho=|\epsilon_n(z)|/\epsilon_n^{(e)}(z)$; 
$z=re^{i\theta}$, $r=10$.
\label{tab:PT.E17}}
\vspace*{0.3cm}
\begin{centering}
  \begin{tabular}{c|ccccccccc}
  \hline
$\theta$ &$\tfrac08\pi$ &$\tfrac18\pi$&$\tfrac28\pi$&$\tfrac38\pi$&
$\tfrac48\pi$&$\tfrac58\pi$&$\tfrac68\pi$&$\tfrac78\pi$&$\tfrac88\pi$   
\\ \hline
   $n=5$ &0.29 &0.30 &0.31 &0.34 &0.13 &0.15 &0.18 &0.25 &0.34 \\
  $n=10$ &0.22 &0.23 &0.24 &0.26 &0.07 &0.09 &0.12 &0.20 &0.37 \\     
  $n=15$ &0.18 &0.18 &0.19 &0.21 &0.05 &0.06 &0.08 &0.12 &0.19  \\ 
\hline
  \end{tabular}\\
\end{centering}
\end{table*}

The upper bound $\chi(n)$ of the variation introduced by Olver is not
very sharp for certain values of the parameters. In his analysis
the estimate of the variation along a certain path follows from
\begin{equation}\label{eq:PT.E17}
\Vp_\P(t^{-n})=\intp \frac{n\,d\tau}{|z+\tau e^{i\phi}|^{n+1} }
=\intp \frac{n\,d\tau}{||z|e^{i(\theta-\phi)}+\tau |^{n+1}},
\end{equation}
where $\theta=\ph{z}$ and $\phi\in[-\frac12\pi,\frac12\pi]$ is 
defined by $\cos\phi=\kappa/|z|$, where $\kappa=|a|$ is introduced in
\pref{eq:PT.E9}. The right-hand side in \pref{eq:PT.E17} 
is estimated by Olver as follows
\begin{equation}\label{eq:PT.E18}
\intp \frac{n\,d\tau}{||z|e^{i(\theta-\phi)}+\tau |^{n+1}}
\le \intp \frac{n\,d\tau}{(|z|^2+\tau^2)^{\frac12n+\frac12}}=
\frac{\chi_(n)}{|z|^n}.
\end{equation}

The right-hand side of \pref{eq:PT.E17} can be written as a
Gauss hypergeometric function, and we find for the variation (along
the same path $\P$)
\begin{equation}
	\begin{array}{ll}\label{eq:PT.E19}
\hbox{${\displaystyle{\Vp_\P(t^{-n})}}$}&
\hbox{${\displaystyle{=\ \frac{n}{|z|^n} 
\intp \frac{d\tau}{(u^2+2\cos(\theta-\phi) u+1)^{\frac12n+\frac12}}}}$}\\
&\\
&
\hbox{${\displaystyle{
=\ \frac{1}{|z|^n} \ \F{\frac12n}{\frac12}{\frac12n+1}{\sin^2(\theta-\phi)}.}}$}\\
	\end{array}
\end{equation}
This result is also given in \cite{Olver:1965:ASW}. The value $\chi(n)$ arises when we replace the argument of the 
$F-$function by unity. With this new value of $\Vp_\P$ we re-computed the 
ratios of Table \pref{tab:PT.E17}, and we give the new ratios in Table 
\pref{tab:PT.E18}. We see that indeed the ratios become larger
in the regions $R_2$ and $R_4$, except when $\theta=\pi$. 
That is,  
the estimates of the remainders become more realistic in these regions.

\begin{table*}
\caption{Ratios as in Table \ref{tab:PT.E17}, now with variations 
$\Vp_\P(t^{-n})$ according to \pref{eq:PT.E19}.
\label{tab:PT.E18}}
\vspace*{0.3cm}
\begin{centering}
  \begin{tabular}{c|ccccccccc}
  \hline
$\theta$ &$\tfrac08\pi$ &$\tfrac18\pi$&$\tfrac28\pi$&$\tfrac38\pi$&
$\tfrac48\pi$&$\tfrac58\pi$&$\tfrac68\pi$&$\tfrac78\pi$&$\tfrac88\pi$   
\\ \hline
   $n=5$ &0.29 &0.30 &0.31 &0.34 &0.38 &0.42 &0.43 &0.41 &0.32 \\
  $n=10$ &0.22 &0.23 &0.24 &0.26 &0.31 &0.35 &0.38 &0.39 &0.33 \\     
  $n=15$ &0.18 &0.18 &0.19 &0.21 &0.25 &0.28 &0.29 &0.27 &0.17  \\ 
\hline
  \end{tabular}\\
\end{centering}
\end{table*}
\section{Uniform expansions in terms of elementary functions}
\label{sec:EA}
We transform the differential equation into a standard form and
distinguish between the cases that there are no real turning points
(as for $a>0$), and that there are two real turning points.
For convenience we consider real parameters. 

\subsection{Positive {\protect\boldmath $a$}}\label{sec:AP}
For $a>0$ no oscillations occur on the real $z-$axis.  
\cite{Olver:1959:UAE} gives expansions that cover
all real $z$. We consider two different modifications, 
one for $z\ge0$ and another one for $z\le0$. These modifications are 
derived in \cite{Temme:2000:NAP}, and we take the same notation 
as in this reference.
  
The function $U(a,z)$ is a solution of the differential equation 
\pref{eq:IN.I1}, and $w(t)=U(\frac12\mu^2,\mu t\sqrt{2})$ satisfies
\begin{equation}
\frac{d^2w}{dt^2} =\mu^4(t^2+1)w.
\end{equation}
The function $W(t)=(t^2+1)^{\frac14} U(\frac12\mu^2,\mu t\sqrt{2})$ 
is a solution of
\begin{equation}\label{eq:UE.T1}
\frac{d^2W}{d\wt\xi^2}=\left[\mu^4+\psi(\wt\xi)\right]W,
\end{equation}
where the relation between $t$ and $\wt\xi$ is given by
\begin{equation}\label{eq:UE.T2}
\wt\xi=\tfrac12t\sqrt{{t^2+1}}+\tfrac12\ln\left[t+\sqrt{{t^2+1}}\right],
\end{equation}
and $\psi(\wt\xi)$ is given by
\begin{equation}\label{eq:UE.T3}
\psi(\wt\xi)=
\frac{2-3t^2}{4(t^2+1)^3}.
\end{equation}
Transformations of this kind are discussed in 
\cite{Olver:1974:ASF}, Chapter 10. The relation in \pref{eq:UE.T2}
follows from the relation
\begin{equation}\label{eq:UE.TD}
\frac{d\wt\xi}{dt}=\sqrt{t^2+1},\quad \wt\xi(0)=0,
\end{equation}
that Olver used in a Liouville-Green transformation.

The quantity $\wt{F}$ in
$W=e^{-\mu^2\wt\xi}\wt{F}$ is a solution of 
\begin{equation}\label{eq:UE.T4}
\frac{d^2\wt{F}}{d\wt\xi^2} - 2\mu^2 \frac{d\wt{F}}{d\wt\xi} -\psi(\wt\xi) \wt{F}=0.
\end{equation}
It is convenient to introduce another parameter, $\wt\tau$, by writing
\begin{equation}\label{eq:UE.T5}
\wt\tau=\tfrac12\left[\frac{t}{\sqrt{{t^2+1}}}-1\right].
\end{equation}
We have
\begin{equation}\label{eq:UE.T6}
\frac{d\wt\xi}{d\wt\tau}=\frac1{8\wt\tau^2(1+\wt\tau)^2}
\end{equation}
and equation \pref{eq:UE.T4} becomes in terms of $\wt\tau$:
\begin{equation}\label{eq:UE.T7}
16\wt\tau^2(\wt\tau+1)^2\frac{d^2\wt{F}}{d\wt\tau^2} +
\left[32\wt\tau(2\wt\tau^2+3\wt\tau+1)-4\mu^2\right]\frac{d\wt{F}}{d\wt\tau}
+(20\wt\tau^2+20\wt\tau+3) \wt{F}=0.
\end{equation}
\subsubsection{\protect\boldmath $z\ge0$}\label{sec:EE.PZ}

We give an asymptotic expansion of the $U-$ function for $a$ large 
and positive, that holds uniformly for $z\ge0$. We write
\begin{equation}\label{eq:UE.T8}
U\left(\tfrac12\mu^2,\mu t\sqrt{2}\right)=
\frac{e^{-\mu^2\wt\xi}}{\sqrt{2}\mu h(\mu) (t^2+1)^{\frac14}}\wt{F}_\mu(\wt\tau),
\end{equation}
where $\wt{F}_\mu(\wt\tau)$ satisfies equation \pref{eq:UE.T7}
and is expanded in the form
\begin{equation}\label{eq:UE.T9}
\wt{F}_\mu(\wt\tau)\sim\sum_{s=0}^\infty(-1)^s\frac{\phi_s(\wt\tau)}{\mu^{2s}},
\end{equation}
where
\begin{equation}\label{eq:UE.T10}
h(\mu)=2^{-\tfrac14\mu^2-\tfrac14}e^{-\tfrac14\mu^2}
\mu^{\tfrac12\mu^2-\tfrac12}=2^{-\frac12} e^{-\frac12a}a^{\frac12a-\frac14}.
\end{equation}
Substituting 
\pref{eq:UE.T9} into \pref{eq:UE.T7} and prescribing 
\begin{equation}\label{eq:UE.T11}
\phi_0(\wt\tau)=1, \quad \phi_s(0)=0, \quad s\ge1,
\end{equation}
we find that the coefficients $\phi_s(\wt\tau)$ 
are polynomials in  $\wt\tau$ of degree $3s$, and are
given by the recursion relation
\begin{equation}\label{eq:UE.T12}
\phi_{s+1}(\tau)=-4\tau^2(\tau+1)^2\frac{d}{d\tau}\phi_s(\tau)
-\tfrac14\int_0^\tau\left(20u^2+20u+3\right)\phi_{s}(u)\,du.
\end{equation}

For deriving this relation observe that \pref{eq:UE.T7}
can be written in the form
\begin{equation}
\mu^2\frac{d\wt{F}}{d\tau}=4\frac{d}{d\tau}
\left[\tau^2(\tau+1)^2\frac{d\wt{F}}{d\tau}\right]
+\tfrac14(20\tau^2+20\tau+3) \wt{F}.
\end{equation}

The term $h(\mu)$ given in \pref{eq:UE.T10} follows 
from  \pref{eq:PT.E1} and from the condition on
$\phi_s(\wt\tau)$ given in \pref{eq:UE.T11}.

The expansion in \pref{eq:UE.T9} corresponds with the expansion (11.10)
given in \cite{Olver:1959:UAE}. In our expansion one of the parameters 
$a$ or $t$ should be large, in Olver's expansion $a$ should be large.
Both Olver's and our expansions hold uniformly for all real $t$, but we 
prefer for negative values of $t$ a slightly different expansion that will 
be given in the next subsection. The expansions also hold in unbounded 
complex domains. For details we refer to \cite{Olver:1959:UAE}. 

The first few coefficients are 
\begin{equation}
\begin{array}{ll}\label{eq:UE.T13}
\phi_0(\tau)&=   1,\morespace\\
\phi_1(\tau)&=   -\Frac{\tau}{12}  (20\tau^2+30\tau+9),\morespace\\
\phi_2(\tau)&=    \Frac{\tau^2}{288}(6160\tau^4+18480\tau^3+19404\tau^2+8028\tau+945),\morespace\\
\phi_3(\tau)&=   -\Frac{\tau^3}{51840}(27227200\tau^6+ 122522400\tau^5 + 220540320\tau^4+\morespace\\
     &   \quad\quad   200166120\tau^3+ 94064328\tau^2+
          20545650\tau+1403325),
\end{array}
\end{equation}

For the derivative we have
\begin{equation}\label{eq:UE.D4}
U'\left(\tfrac12\mu^2,\mu t\sqrt{2}\right)=
-\frac{(1+t^2)^{\tfrac14}\,e^{-\mu^2{\wt\xi}}}{2h(\mu)}
{\wt G}_a(z),\quad
{\wt G}_a(z)\sim
\sum_{s=0}^\infty\,(-1)^s\,\frac{\psi_s(\wt\tau)}{\mu^{2s}},
\end{equation} 
where the coefficients $\psi_s(\tau)$ can be obtained by formal 
differentiating  \pref{eq:UE.T8} and  \pref{eq:UE.T9}. 
It follows that
\begin{equation}\label{eq:UE.T17}
\psi_s(\tau)=\phi_s(\tau)+2\tau(\tau+1)(2\tau+1)
\phi_{s-1}(\tau)+8\tau^2(\tau+1)^2\frac{d\phi_{s-1}(\tau)}{d\tau},
\end{equation}
$s=0,1,2,\ldots\,$. The first few coefficients are
\begin{equation}
\begin{array}{ll}\label{eq:UE.T18}
\psi_0(\tau)&=   1,\morespace\\
\psi_1(\tau)&=   \Frac{\tau}{12}  (28\tau^2+42\tau+15),\\
\psi_2(\tau)&=   -\Frac{\tau^2}{288}(7280\tau^4+21840\tau^3+23028\tau^2+9684\tau+1215),\morespace\\
\psi_3(\tau)&=\Frac{\tau^3}{51840}(30430400\tau^6+136936800\tau^5+246708000\tau^4+\morespace\\
      &   \quad\quad
224494200\tau^3+106122312\tau^2+23489190\tau+1658475).
\end{array}
\end{equation}
\subsubsection{{\protect\boldmath $z\le0$}}\label{sec:EE.NZ}

For negative $z$ we have
\begin{equation}
\begin{array}{ll}\label{eq:UE.D5}
U\left(\tfrac12\mu^2,-\mu t\sqrt{2}\right)
&=\Frac{\sqrt{{2\pi}}}{\Gamma(\tfrac12+\tfrac12\mu^2)}
\Frac{h(\mu)e^{\mu^2{\wt\xi}}}{(1+t^2)^{\frac14}}\,
\wt P_\mu(t),\\ 
U'\left(\tfrac12\mu^2,-\mu t\sqrt{2}\right)
&=-\Frac{\sqrt{\pi}\,\mu h(\mu)}{\Gamma(\tfrac12+\tfrac12\mu^2)}
 e^{\mu^2{\wt\xi}}(1+t^2)^{\tfrac14}\,
\wt Q_\mu(t),
\end{array}
\end{equation}
where
\begin{equation}\label{eq:UE.D6}
\wt P_\mu(t)\sim\sum_{s=0}^\infty \frac{\phi_s(\wt\tau)}{\mu^{2s}},\quad
\wt Q_\mu(t)\sim\sum_{s=0}^\infty \frac{\psi_s(\wt\tau)}{\mu^{2s}}.
\end{equation}
Again, these expansions are valid when one or both parameters $a$ and $t$ 
are large.

The functions $\wt F_\mu(t), \wt G_\mu(t), 
\wt P_\mu(t)$ and $\wt Q_\mu(t)$ 
satisfy the following exact relation:
\begin{equation}\label{eq:UE.D7}
\wt F_\mu(t)\wt Q_\mu(t)+
\wt G_\mu(t)\wt P_\mu(t)=2.
\end{equation}

The relation in \pref{eq:IN.I17} can be used for obtaining expansions 
for $V(a,z)$ and its derivative.
\subsubsection{Error bounds of the expansions}\label{sec:EE.RD}
We apply Theorem 3.1 of \cite{Olver:1974:ASF}, page 366, and write the 
expansion in \pref{eq:UE.T9} with a remainder. 
For $n=0,1,2,\ldots$ we have
\begin{equation} \label{eq:UE.D8}
\wt F_\mu(t)=\sum_{s=0}^{n-1}(-1)^s 
\frac{\phi_s(\wt\tau)}{\mu^{2s}}
+\wt R_{n}(\mu,t).
\end{equation}
The remainder $\wt R_{n}(\mu,t)$ can be bounded as follows
\begin{equation}\label{eq:UE.D9}
|\wt R_{n}(\mu,t)|\le  
\exp\left[
\frac{2\Vp_{\wt\xi, \infty}(\phi_1)}{\mu^2}\right]\, 
\frac{\Vp_{\wt\xi, \infty}(\phi_n)}{\mu^{2n}}, 
\end{equation}
where we take into account that we consider positive $t$ and $\mu$
(see also Exercise 3.1 on page 367 of \cite{Olver:1974:ASF}).

We need the variation of the coefficients $\phi_s$ for $t\ge0$
which corresponds with  $\wt\xi\ge0$ and  
$-\frac12\le\wt\tau\le0$; cf. \pref{eq:UE.T5}. We have
\begin{equation}\label{eq:UE.D10}
\Vp_{\wt\xi, \infty}(\phi_n)=\int_{\wt 
\tau}^{0}|\phi_n'(\tau)|\,d\tau. 
\end{equation}
For negative argument the coefficients $\phi_s$ oscillate; 
see Figure \ref{fig:fig2}.

\begin{center}
\epsfxsize=8cm \epsfbox{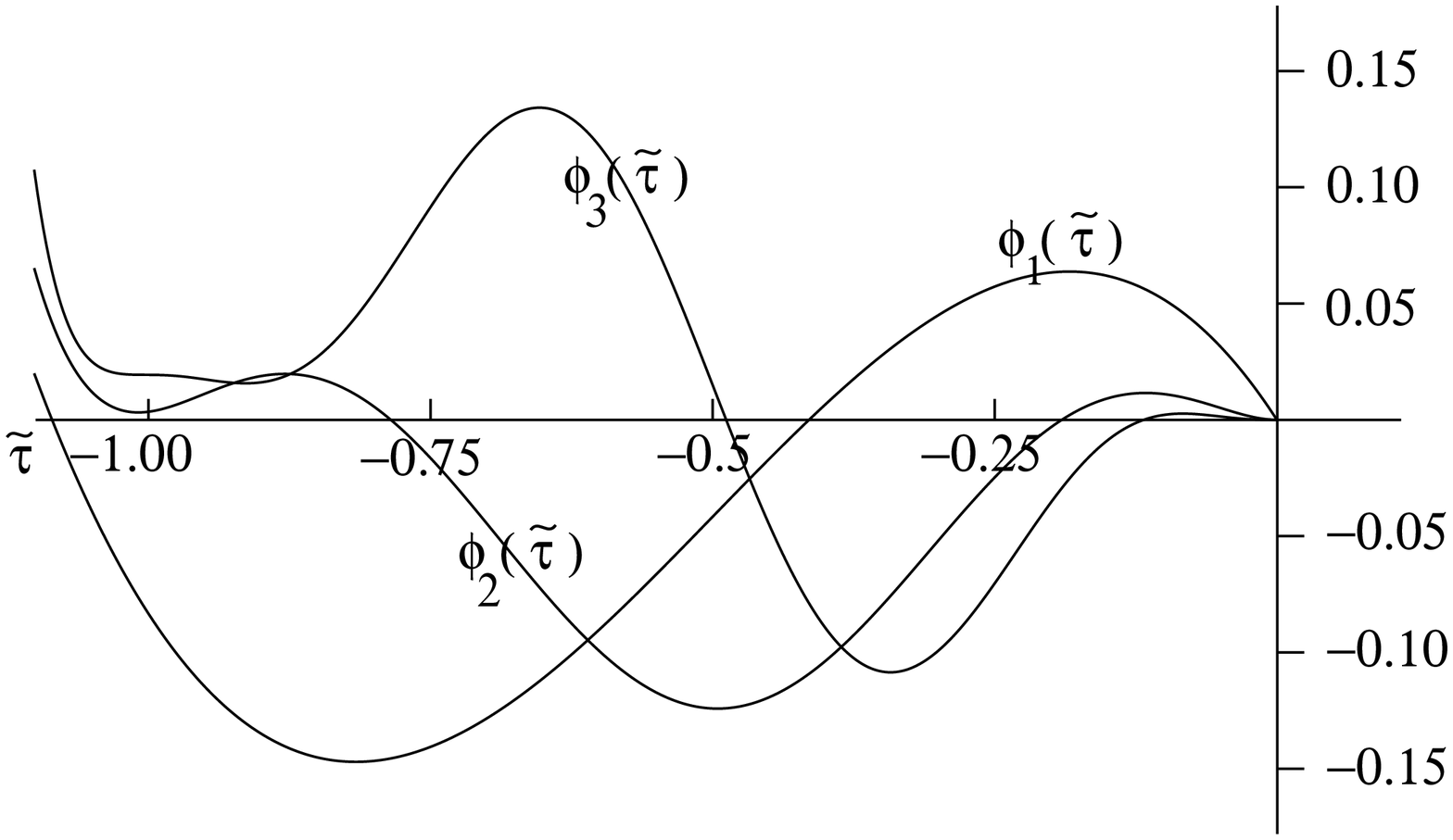}\label{fig:fig2}
\vspace*{0.5cm}

{\bf Figure \ref{fig:fig2}.}\ Graphs of $\phi_s(\wt\tau)$, $s=1,2,3$, $\wt\tau
\in[-1,0]$.
\end{center}

In Table \ref{tab:UE.Tpp} we show the ratios 
$|R_{n}(\mu,t)|/R_{n}^{(e)}(\mu,t)$ (with $n=3$), where $R_{n}^{(e)}(\mu,t)$ 
is the right-hand side of \pref{eq:UE.D9}. We see that the estimates
for small values of $t$ are much too large. An explanation is that 
for small $t$ (that is $\wt\tau$ close to $-\frac12$), the variations 
in \pref{eq:UE.D10} are calculated over a larger interval than when $t$ is 
large.

\begin{table*}
\caption{Ratios $|R_{n}(\mu,t)|/R_{n}^{(e)}(\mu,t)$; $z=2t\sqrt{a}\,, n=3$.
\label{tab:UE.Tpp}}
\vspace*{0.3cm}
\begin{centering}
  \begin{tabular}{c|ccccccccc}
  \hline
$t$ &$0.0$ &$1.0$&$2.5$&$5.0$&
$10$&$25$&$50$   
\\ \hline
$a= 1$   &.21493  &.14455  &.84677  &.94186  &.98360  &.99728  &.99932\\
$a= 5$   &.06142  &.43256  &.96494  &.98773  &.99667  &.99945  &.99986\\
$a= 10$  &.00343  &.50123  &.98214  &.99382  &.99833  &.99973  &.99993\\
$a= 50$  &.04921  &.56597  &.99637  &.99876  &.99967  &.99995  &.99999\\
$a= 100$ &.05601  &.57478  &.99818  &.99938  &.99983  &.99997  &.99999\\
\hline
  \end{tabular}\\
\end{centering}
\end{table*}

For $z\le0$ we consider the remainder in the expansion of 
\pref{eq:UE.D6} and we write
\begin{equation}\label{eq:UE.11}
\wt P_\mu(t)=\sum_{s=0}^{n-1} \frac{\phi_s(\wt\tau)}{\mu^{2s}}+
\wt R_{n}(\mu,t).
\end{equation}
In the present case we have the bound 
\begin{equation}\label{eq:UE.D12}
|\wt R_{n}(\mu,t)|\le  
\exp\left[
\frac{2\Vp_{-\infty,\wt\xi}(\phi_1)}{\mu^2}\right]\, 
\frac{\Vp_{-\infty,\wt\xi}(\phi_n)}{\mu^{2n}}, 
\end{equation}
where
\begin{equation}\label{eq:UE.D13}
\Vp_{-\infty,\wt\xi}(\phi_n)=\int_{-1}^{\wt \tau}|\phi_n'(\tau)|\,d\tau. 
\end{equation}

In Table \ref{tab:UE.Tpn} we again show the ratios 
$|R_{n}(\mu,t)|/R_{n}^{(e)}(\mu,t)$ (with $n=3$), 
where in the present case $R_{n}^{(e)}(\mu,t)$ 
is the right-hand side of \pref{eq:UE.D12}. We see that the ratios
in general are smaller than when $z\ge0$. An explanation is that 
the variations include the contributions from the interval 
$[-1,-\frac12]$. The parameter $\wt\tau$ remains in the interval 
$[-\frac12,0]$, however. 

\begin{table*}
\caption{Ratios $|R_{n}(\mu,t)|/R_{n}^{(e)}(\mu,t)$; $z=-2t\sqrt{a}\,, n=3$.
\label{tab:UE.Tpn}}
\vspace*{0.3cm}
\begin{centering}
  \begin{tabular}{c|ccccccccc}
  \hline
$t$ &$0.0$ &$1.0$&$2.5$&$5.0$&
$10$&$25$&$50$   
\\ \hline
$a= 1$   &.29041   &.04469  &.87352  &.76079  &.72513  &.71493  &.71347\\
$a= 5$   &.17780   &.02071  &.96996  &.94637  &.93771  &.93509  &.93471\\
$a= 10$  &.12433   &.01817  &.98467  &.97279  &.96835  &.96700  &.96680\\
$a= 50$  &.07476   &.01644  &.99689  &.99449  &.99359  &.99331  &.99327\\
$a= 100$ &.06829   &.01624  &.99844  &.99724  &.99679  &.99665  &.99663\\
\hline
  \end{tabular}\\
\end{centering}
\end{table*}
\subsubsection{Upper bounds for the variations of
{\protect\boldmath $\phi_s(\wt\tau)$}}\label{sec:EE.UB}
The variations of the coefficients $\phi_s(\wt\tau)$ used in \pref{eq:UE.D10}
and \pref{eq:UE.D13} can be computed by numerical quadrature of the integrals,
but for real values it is convenient to use the zeros of the polynomials
$\phi_s'(\wt\tau)$. For example,  $\phi_1'(\wt\tau)$ has zeros at
$t_1 = -0.816$ and $t_2 = -0.184$. Hence, for 
$\wt\tau\in[-\frac12,0]$ the variation in \pref{eq:UE.D10} follows from
\begin{equation}\label{eq:UE.N10}
\Vp_{\wt\xi, \infty}(\phi_1)=\int_{\wt 
\tau}^{0}|\phi_1'(\tau)|\,d\tau = \left\{
\begin{array}{ll}
\phi_1(\wt\tau)              &\mbox{if\  } \wt\tau\in[t_2,0],\\ 
2\phi_1(t_2)-\phi_1(\wt\tau) &\mbox{if\  } \wt\tau\in[-\frac12,t_2].
\end{array}
\right.
\end{equation}
The computation of the zeros of $\phi_s'(\wt\tau)$, however, may be not
efficient in an algorithm. We can avoid this by constructing 
for the variations upper bounds as functions of $\wt\tau$.
For example, we found for the first few $\phi_s(\wt\tau)$ the following 
simple upper  bounds for $\wt\tau\in[-\frac12,0]$:
 \begin{eqnarray}\label{eq:UE.N11}
\Vp_{\wt\xi, \infty}(\phi_1) &\le& 
\Frac{-3\,\wt\tau}{4(1+4.8\wt\tau^2)},\nonumber\\
\Vp_{\wt\xi, \infty}(\phi_2) &\le& 
\Frac{105\,\wt\tau^2}{32(1+18\wt\tau^2)},\\
\Vp_{\wt\xi, \infty}(\phi_3) &\le& 
\Frac{-3465\,\wt\tau^3}{128(1+52\wt\tau^2)}.\nonumber
\end{eqnarray}
The bounds fit at the origin, and are slightly larger at $\wt=-\frac12$. 

In Figure \ref{fig:fig3} we give the graphs of $\phi_s(\wt\tau)$, $s=1$
(left), $s=2$ (right), and $s=3$ for $\wt\tau\in[-\frac12,0]$.  The first
graph from the bottom is for $\phi_s(\wt\tau)$, the second one is for the
variation in \pref{eq:UE.D10}, the third one is for the upper bound given
in \pref{eq:UE.N11}.  The forth one is for the variation in
\pref{eq:UE.D13}, and the fifth one is the upper bound.  Those two are just
equal to the second and third plus the variation over $[-1,-\frac12]$,
which equals $0.1692$, $0.1602$, $0.2415$, for $s=1,2,3$, respectively.
Observe that in  \pref{eq:UE.D13} we consider $\wt\tau\in[-\frac12,0]$,
and that the variation contains contributions from $[-1,-\frac12]$.

\begin{center}
\epsfxsize=12cm \epsfbox{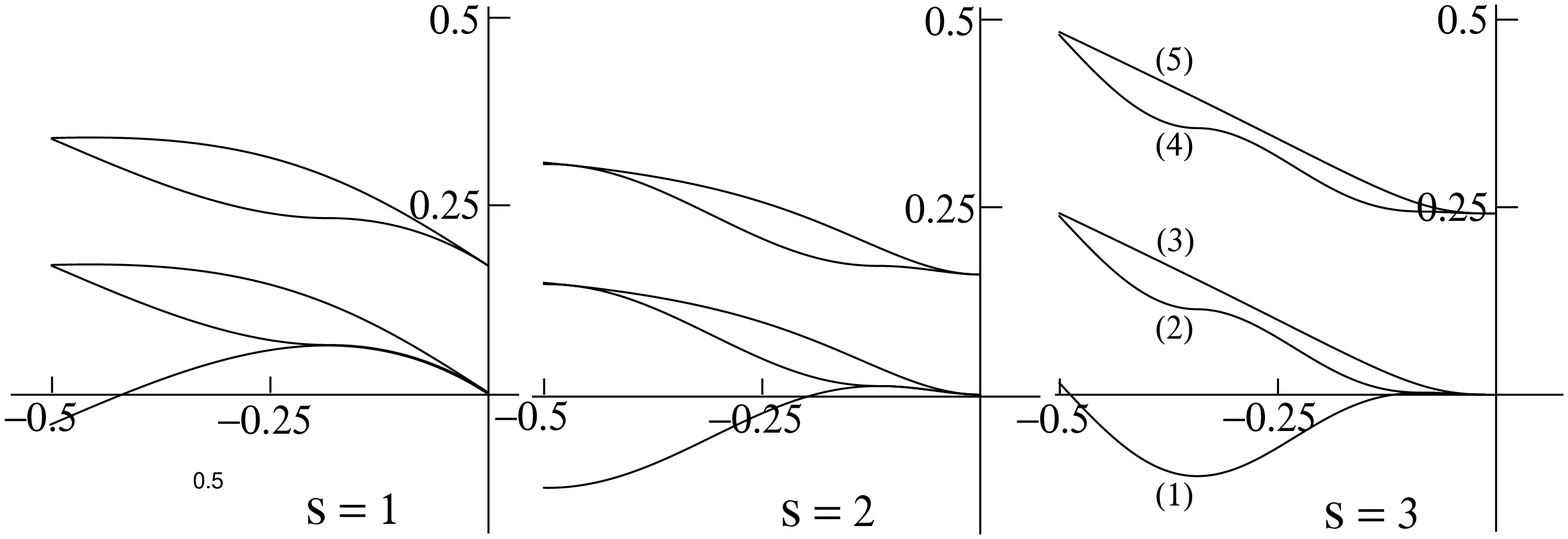}\label{fig:fig3}
\vspace*{0.5cm}

{\bf Figure \ref{fig:fig3}.}\ Graphs of $\phi_s(\wt\tau)$, $s=1,2,3$, $\wt\tau
\in[-\frac12,0]$, variations and upper bounds; see the text for more details.
\end{center}

\subsection{Negative {\protect\boldmath $a$}}\label{sec:AN}
In this case we consider the function $U(-\frac12\mu^2,\mu t\sqrt{2})$. 
This function satisfies the differential equation
\begin{equation}
\frac{d^2w}{dt^2} =\mu^4(t^2-1)w.
\end{equation}
For expansions in terms of elementary functions three intervals should be 
distinguished: $(-\infty,-1-\delta]$, $[-1+\delta,1-\delta]$ and
$[1+\delta,\infty)$ (the turning points $t=\pm 1$ should be avoided).
We only consider the interval $[1+\delta,\infty)$; for the other 
intervals we refer to \cite{Temme:2000:NAP}.

\subsubsection{{\protect\boldmath $t>1$.}}\label{sec:EE.T1}

We write (see \cite{Temme:2000:NAP})
\begin{equation}\label{eq:UE.T8N}
U\left(-\tfrac12\mu^2,\mu t\sqrt{2}\right)=
\frac{h(\mu)\,e^{-\mu^2\xi}}{(t^2-1)^{\frac14}}F_\mu(\tau),
\end{equation}
where $F_\mu(\tau)$ is expanded in the form
\begin{equation}\label{eq:UE.T9N}
F_\mu(\tau)\sim\sum_{s=0}^\infty\frac{\phi_s(\tau)}{\mu^{2s}},
\end{equation}
and where $h(\mu)$ is defined in \pref{eq:UE.T10}, 
\begin{equation}\label{eq:UE.T5N}
\tau=\tfrac12\left[\frac{t}{\sqrt{{t^2-1}}}-1\right],
\end{equation}
\begin{equation}\label{eq:UE.T2N}
\xi=\tfrac12t\sqrt{{t^2-1}}-\tfrac12\ln\left[t+\sqrt{{t^2-1}}\right],
\end{equation}
and the coefficients $\phi_s(\tau)$ are as in \pref{eq:UE.T9}; see
also \pref{eq:UE.T12} and \pref{eq:UE.T13}.

The analysis for deriving \pref{eq:UE.T9} is similar to that of the case 
$a>0$. The function $F(\tau)$ satisfies the equation (see \pref{eq:UE.T7})
\begin{equation}
\mu^2\frac{dF}{d\tau}=
-4\frac{d}{d\tau}\left[\tau^2(\tau+1)^2\frac{dF}{d\tau}\right]
-\tfrac14(20\tau^2+20\tau+3) F.
\end{equation}

For the function $V(a,z)$ we have 
\begin{equation}\label{eq:UE.T14}
V\left(-\tfrac12\mu^2,\mu t\sqrt{2}\right)=
 \frac{       e^{\mu^2\xi}}{\mu\,\sqrt{{\pi}}\,h(\mu)(t^2-1)^{\frac14}}
\,P_\mu(t),\quad
P_\mu(t)\sim\sum_{s=0}^\infty(-1)^s\frac{\phi_s(\tau)}{\mu^{2s}},
\end{equation}
where the $\phi_s(\tau)$ are the same as in \pref{eq:UE.T9}.

For the derivatives we have
\begin{equation}\label{eq:UE.T15}
U'\left(-\tfrac12\mu^2,\mu t\sqrt{2}\right)=
-\frac{\mu}{\sqrt{2}}h(\mu)(t^2-1)^{\tfrac14}e^{-\mu^2\xi}\,G_\mu(t),
\quad G_\mu(t)\sim
\sum_{s=0}^\infty\frac{\psi_s(\tau)}{\mu^{2s}} 
\end{equation}
and
\begin{equation}\label{eq:UE.T16}
V'\left(-\tfrac12\mu^2,\mu t\sqrt{2}\right)=
\frac{(t^2-1)^{\tfrac14}e^{\mu^2\xi}}{\sqrt{{2\pi}} h(\mu)}\,Q_\mu(t),
\quad Q_\mu(t)\sim
\sum_{s=0}^\infty(-1)^s\frac{\psi_s(\tau)}{\mu^{2s}}. 
\end{equation}

The coefficients $\psi_s$ are the same as in \pref{eq:UE.D4}; see also
\pref{eq:UE.T17} and \pref{eq:UE.T18}.

\subsubsection{Error bounds of the expansion}\label{sec:EA.RA}
We write the expansion in \pref{eq:UE.T9} in the form
\begin{equation}\label{eq:UE.T20}
F_\mu(t)={\displaystyle{\sum_{s=0}^{n-1}\frac{\phi_s(\tau)}{\mu^{2s}}}} 
+R_{n}(\mu,t).
\end{equation}

For the present values of the parameters the remainder $R_{n}(\mu,t)$ 
can be bounded as follows,
see \cite{Olver:1974:ASF} (page 366),
\begin{equation}\label{eq:UE.T21}
|R_{n}(\mu,t)|\le  
\exp\left[\frac{2\Vp_{\infty,\xi}(\phi_1)}{|\mu^2|}\right]\,
\frac{\Vp_{\infty,\xi}(\phi_n)}{|\mu^2|^n}.
\end{equation}

In Table \ref{tab:UE.T22} we show the ratios 
$|R_{n}(\mu,t)|/R_{n}^{(e)}(\mu,t)$ (with $n=3$), where $R_{n}^{(e)}(\mu,t)$ 
is the right-hand side of \pref{eq:UE.T21}.

\begin{table*}
\caption{Ratios $|R_{n}(\mu,t)|/R_{n}^{(e)}(\mu,t)$; $z=2t\sqrt{a}\,, n=3$.
\label{tab:UE.T22}}
\vspace*{0.3cm}
\begin{centering}
  \begin{tabular}{c|ccccccccc}
  \hline
$t$ &$1.5$ &$2.0$&$3.0$&$5.0$&
$10$&$20$&$50$   
\\ \hline
$a=- 1$   &.29990& .57546& .80676& .93078& .98282& .99572& .99932\\
$a=- 5$   &.69226& .86898& .95344& .98522& .99652& .99914& .99986\\
$a=- 10$  &.81624& .92930& .97608& .99256& .99826& .99956& .99994\\
$a=- 50$  &.95602& .98488& .99510& .99850& .99964& .99992& .99998\\
$a=- 100$ &.97744& .99236& .99754& .99924& .99982& .99996& 1.0000\\
\hline
  \end{tabular}\\
\end{centering}
\end{table*}

From the first few $\phi_s(\tau)$ given in \pref{eq:UE.T13}
and from the recursion relation in \pref{eq:UE.T12} it follows that all
coefficients in these polynomials have the sign of $(-1)^s$. Hence, for 
$t>1$, that is, $\tau \ge 0$, the variations in
\pref{eq:UE.T21} can be easily obtained. We have for $n\ge1$, using 
\pref{eq:PT.E11},
\begin{equation}\label{eq:UE.T26} 
\Vp_{\infty,\xi}(\phi_n)=(-1)^s\int_0^\tau\phi_n'(\sigma)\,d\sigma=
(-1)^s\phi_n(\tau)=|\phi_n(\tau)|.
\end{equation}

\section{Error bounds by using integrals}\label{sec:EE.RI}
The construction of error bounds of remainders in uniform asymptotic
expansions is available now by Olver's work on differential equations.  In
this section we consider a method for expansion \pref{eq:UE.T9} by using
an integral representation of the function $U(a,z)$.  For this approach it is
convenient to concentrate on large positive values of $z$, and to construct
an expansion that holds uniformly with respect to $a\in [0,\infty)$.
This expansion reduces to the Poincar\'e-type expansion in \pref{eq:PT.E1} when
$a$ is fixed after expanding the quantities in the expansion 
that depend on $\la=a/z^2$ for small 
values of this parameter. In fact, by writing $\mu=z\sqrt{2\la}$ and 
$t=1/(2\sqrt{\la})$ the same can be done in the expansion 
given by \pref{eq:UE.T8} and \pref{eq:UE.T9}.

\subsection{An integration by parts procedure}
\label{sec:EE.IR}
We summarize from \cite{Temme:2000:NAP} and start with the integral representation
\begin{equation}
U(a,z)=\frac{e^{-\tfrac14z^2}}{\Gamma(a+\tfrac12)}
\intp w^{a} e^{-\tfrac12w^2-zw}\,\frac{dw}{\sqrt{w}},\quad a>-\tfrac12\label{i1}
\end{equation}
which we write in the form
\begin{equation}
U(a,z)=\frac{z^{a+\tfrac12}\,e^{-\tfrac14z^2}}{\Gamma(a+\tfrac12)}
\intp e^{-z^2\phi(w)}\, \frac{dw}{\sqrt{w}},\label{i2}
\end{equation}
where 
\begin{equation}
\phi(w)=w+\tfrac12w^2-\lambda\ln w, \quad 
\lambda=\frac{a}{z^2}=\frac{1}{4t^2},\label{i50}
\end{equation}
where $t$ is used earlier in the notation $U(\frac12\mu^2,\mu t\sqrt{2})$. 
The positive saddle point $w_0$ of $\phi(w)$ is
\begin{equation}
 w_0=\tfrac12\left[\sqrt{{1+4\lambda}}-1\right].\label{i5}
\end{equation}

A standard form of \pref{i2} is obtained by using the transformation
\begin{equation}
\phi(w)=s-\lambda\ln s+A,\label{i6}
\end{equation}
where $A$ does not depend on $s$ or $w$; we prescribe that $w=0$
should correspond with $s=0$, and $w=w_0$  with $s=\lambda$, the saddle point
in the $s-$plane. This gives
\begin{equation}
A=\tfrac12w_0^2+w_0-\lambda\ln w_0-\lambda+\lambda\ln\lambda,\label{i9}
\end{equation}
\begin{equation}
U(a,z)=\frac{z^{a+\tfrac12}\,e^{-\tfrac14z^2-Az^2}}{(1+4\lambda)^{\frac14}
\Gamma(a+\tfrac12)}\intp
s^{a}e^{-z^2s} f(s)\,\frac{ds}{\sqrt{s}},\label{i7}
\end{equation}
where
\begin{equation}
f(s)=(1+4\lambda)^{\tfrac14}\sqrt{{\frac sw }}\,\frac{dw}{ds}
=(1+4\lambda)^{\tfrac14}
\sqrt{{\frac ws }}\,\frac{s-\lambda}{w^2+w-\lambda}.\label{i8}
\end{equation}
By normalizing with the quantity $(1+4\lambda)^{\tfrac14}$ we obtain $f(\lambda)=1$,
as can be verified from \pref{i8} and a limiting process (using
l'H\^opital's rule). 

For $\la\to0$ the saddle point $w_0$ tends to zero, and the mapping becomes
\begin{equation}\label{ii13}
\tfrac12w^2+w  =s.
\end{equation}
It is not difficult to verify that  for $\la=0$ we have
\begin{equation}\label{ii14}
f(s)=\sqrt{\frac{1+\sqrt{1+2s}}{2(1+2s)}}.
\end{equation}
If $\la\ne0$ the transformation \pref{i6} can also be written 
\begin{equation}\label{ii15}
\la\,w= w_0\,s\,e^{\frac{1}{\la}(\frac12w^2+w-s-\frac12w_0^2-w_0+\la)},
\end{equation}
in which form no logarithms occur.

We introduce a sequence of functions $\{f_k\}$ with
$f_0(s)=f(s)$ and 
\begin{equation}
f_{k+1}(s)=\sqrt{{s}}\,\frac{d\ }{ds}
\left[\sqrt{{s}}\frac{f_k(s)-f_k(\lambda)}{s-\lambda}\right],\quad
k=0,1,2,\ldots\ .\label{i14}
\end{equation}

The expansion in \pref{eq:UE.T9} can be obtained by using an integration by
parts procedure. 
Consider the integral
\begin{equation}
F_a(z)=\frac1{\Gamma(a+\tfrac12)}
\intp s^{a}e^{-z^2s} f(s)\,\frac{ds}{\sqrt{s}},\label{i21}
\end{equation}
We have (with $\lambda=a/z^2$)
\begin{eqnarray*}
F_a(z)
&=& z^{-2a-1}f(\lambda) +
\Frac{1}{\Gamma(a+\tfrac12)}
\intp s^{a}e^{-z^2s}[f(s)-f(\lambda)] \,\frac{ds}{\sqrt{s}}\\
&=& z^{-2a-1}f(\lambda)  -
\Frac{1}{z^2\Gamma(a+\tfrac12)}
\intp \sqrt{s} \Frac{f(s)-f(\lambda)}{s-\lambda} 
\,de^{-z^2(s-\lambda\ln s)}\\
&=& z^{-2a-1}f(\lambda) +  \Frac1{z^2\Gamma(a+\tfrac12)}
\intp s^{a}e^{-z^2s} f_1(s)\,\frac{ds}{\sqrt{s}},
\end{eqnarray*}
where $f_1$ is given in \pref{i14} with $f_0=f$. Repeating this
procedure we obtain 
\begin{equation}
U(a,z)\sim 
\frac{e^{-\tfrac14z^2-Az^2}}{z^{a+\tfrac12}\,(1+4\lambda)^{\frac14}}
\sum_{k=0}^\infty \frac{f_k(\lambda)}{z^{2k}}.\label{i13}
\end{equation}
The factors in front of the series in \pref{i13} and \pref{eq:UE.T9}
are the same. This can be verified by using $a=\tfrac12\mu^2$
and $z=\mu\sqrt{2} t$.
Also, the series correspond termwise with each other,
the relation between the coefficients being
\begin{equation}\label{i144}
\phi_k(\wt\tau)= (-1)^k\,(2\lambda)^{k}\,f_k(\lambda),\quad
\wt\tau=\tfrac12\left[\frac{1}{\sqrt{{4\lambda+1}}}-1\right].
\end{equation}
For example, we have
\begin{equation}
f_0(\lambda) =   1,\quad f_1(\lambda)=
-\frac{(2\wt\tau+1)^2}{24(\wt\tau+1)}  
(20\wt\tau^2+30\wt\tau+9).
\label{i32}
\end{equation}

We write \pref{i13} with a remainder:
\begin{equation}
U(a,z)= 
\frac{e^{-\tfrac14z^2-Az^2}}{z^{a+\tfrac12}\,(1+4\lambda)^{\frac14}}
\left[\sum_{k=0}^{n-1} \frac{f_k(\lambda)}{z^{2k}}+
\frac1{z^{2n}}R_n(a,z)\right],\label{i15}
\end{equation}
where
\begin{equation}
R_n(a,z)=\frac{z^{2a+1}}{\Gamma(a+\frac12)}
\intp s^{a}e^{-z^2s} f_n(s)\,\frac{ds}{\sqrt{s}}.\label{i16}
\end{equation}

\subsection{Bounding the remainder}\label{sec:EE:BR}
From \pref{i6} and \pref{i8} we infer that $f(s)=\bo(s^{-1/4})$ 
as $s\to\infty$ (see also \pref{ii14} for a simple verification when 
$\la=0$). So, $f(s)$ is bounded on $[0,\infty)$.
Further, we can prove (which will not be done here) from \pref{i14} 
that $f_n(s)=\bo(s^{-1/4})$, $n\ge1$. Hence, the functions $f_n(s)$ can be bounded
$|f_n(s)|\le M_n(\lambda)$, for $s\ge0$. Using such a bound in \pref{i16}
will indeed give an upper bound for 
$|R_n(a,z)|$, but this bound may be not realistic.

A much better bound will be obtained by estimating $|f_n(s)|$ accurately
in a small interval around $s=\la$, and accepting a rough estimate for
other $s-$values. This can be achieved by  using a "weight function" 
$w_n(s,\la)$, and by writing
\begin{equation}
|f_n(s)|\le \left[|f_n(\la)|+M_n(\la)\right]\, w_n(s,\la),\label{i17}
\end{equation}
where, for example, we take
\begin{equation}
 w_n(s,\la)=\left[(s/\lambda)^{-\lambda} 
e^{s-\lambda}\right]^{\sigma_n}.\label{i17a}
\end{equation}
We have $w_n(\la,\la)=1$, that is, for the $s-$value where the dominant part 
$s^{a}e^{-z^2s}$ of 
the integrand in \pref{i16} assumes its maximal value
when $a$ and $z$ are large. We try to find  $M_n(\la)>0$ and 
$\sigma_n\ge0$ such that \pref{i17} holds for all $s\ge0$.
Then we obtain the bound
\begin{equation}
|R_n(a,z)|\le  \left[|f_n(\la)|+M_n(\la)\right]\,  S_n(a,z),\label{i18}
\end{equation}
where
\begin{equation}
S_n(a,z) = a^{\lambda \sigma_n} e^{-\lambda \sigma_n}
\left(1-\frac{\sigma_n}{z^2}\right)^{\lambda\sigma_n-a-\frac12}
\frac{\Gamma(a+\frac12-\lambda\sigma_n)}{\Gamma(a+\frac12)}.\label{i19}
\end{equation}
For $S_n(a,z)$ we need the conditions 
\begin{equation}
 z^2>\sigma_n,\quad a+\tfrac12>\lambda\sigma_n.
\end{equation}
The quantity $S_n(a,z)$ is close to unity when $a + z$ is large; because 
$f_n$ is bounded the value of $\sigma_n$ will be small. Numerical 
calculations show that for $\sigma_n=1$ and $z\ge3, a\ge1$ the maximal 
value of $S_n(a,z)$ is smaller than $1.062$\,.

In Figure \ref{fig:fig4} we show the graph of 
$(1+5\la)\left[|f_n(\la)|+M_n(\la)\right]\,$ (the factor 
$(1+5\la)$ is chosen because of scaling)
when we take $\sigma_1=1$ in 
\pref{i17}. We have $M_1(0)=0$, $f_1(0)=-3/8$  and 
$f_1(\la)\sim 1/(48\la)$ for large $\la$. We also draw the graph of 
\begin{equation} \label{rho1}
\rho_1(\la)=  \frac{ |R_1(a,z)|}{\left[|f_1(\la)|+M_1(\la)\right]\, S_1(a,z)},
\end{equation}
the ratio of the exact error and the estimated error. We see a sharp dip 
at $\la=8.3176\ldots$, which is a zero of $f_1(\la)$; for this value of 
$\la$ the
asymptotic approximation improves, as expected. For large $\la$ the 
quantity $\rho_1$ tends to 1.

We computed $M_1$ in  $\left[|f_n(\la)|+M_n(\la)\right]$ as the 
$\inf|f_1(s)/w_1(s,\la)|$, $s\ge0$ and this gives a continuous function 
$M_1(\la)$, which need not to be smooth.
For example there is a noticeable non-smooth behavior near $\la=11$, 
because $\inf|f_1(s)/w_1(s,\la)|$ occurs for $\la<11$ in a different
$s-$domain than for $\la>11$.

\vspace*{0.5cm}

\begin{center}
\epsfxsize=8cm \epsfbox{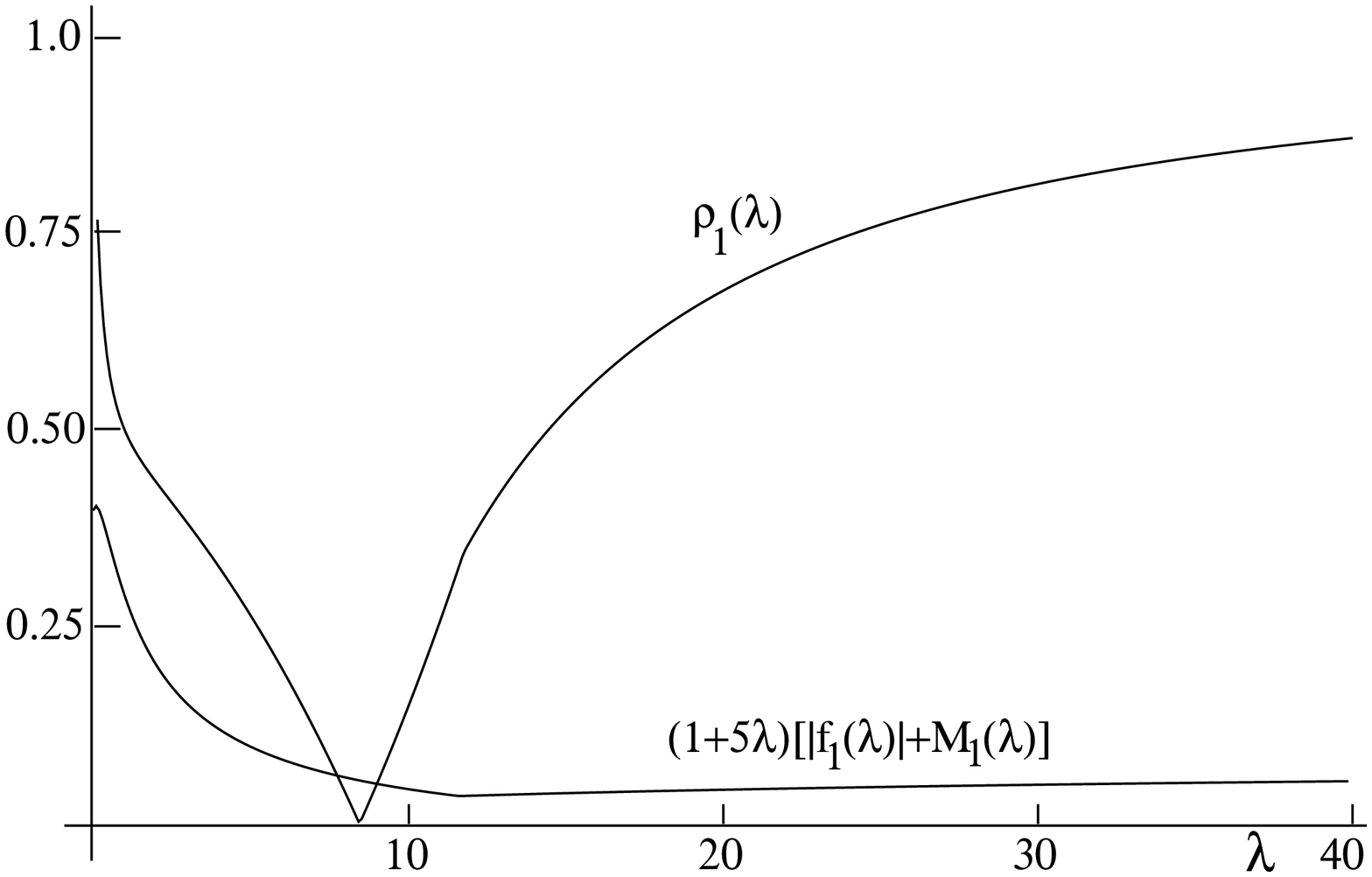}\label{fig:fig4}
\vspace*{0.5cm}

{\bf Figure \ref{fig:fig4}.}\ 
The graph of $(1+5\la)[|f_1(\la)|+M_1(\la)]$ with
$M_1(\la)$ introduced  in \pref{i17}
(with  $\sigma_1=1$), and $\rho_1(\la)$ from \pref{rho1}.
\end{center}

\subsection{Bounding the remainder by using Cauchy-type integrals
}\label{sec:EE:CT}
Computing the functions $f_n(s)$ by using formula  \pref{i14} is quite 
difficult, also when we use computer algebra. The representations contain 
derivatives and removable singularities at $s=\la$. In particular the 
poles are very inconvenient when computing the functions $f_n(s)$ near 
$s=\la$.

It is possible, however, to represent $f_n(s)$ as a Cauchy-type integral.
The mapping in \pref{i6} is singular at $w=w_-=-1-w_0$, the negative saddle point  
of $\phi(w)$ defined in \pref{i50}. The corresponding $s-$value $s_-$ is negative.
If $\la=0$ then $w_-=-1$ and the 
corresponding $s-$value is $s_=-\frac12$. See \pref{ii14}, where indeed 
$f(s)$ shows a singularity at this point. For large values of $\la$ we 
have the estimate (see \cite{Temme:2000:NAP}, formula (4.45))
\begin{equation} \label{sm}
s_-=s(w_-)\sim-\la\left[0.2785+\frac{0.4356}{\sqrt{\la}}\right].
\end{equation}

For constructing the Cauchy-type integrals we use the property that 
the functions $f_n(s)$ are
analytic functions in a domain $\D$ in the half plane $\Re s> s_-$; in particular,
$\D$ contains  a neighborhood of the positive real axis.

As in \cite{OldeDaalhuis:1994:UAT}, we start from
\begin{equation} \label{q0}
f_n(s)=\frac{1}{2\pi i} \int_{\C} Q_0(\sigma,\la,s) f_n(\sigma)\,d\sigma,
\quad Q_0(\sigma,\la,s) =\frac{1}{\sigma-s},
\end{equation}
where $s\in\D$ and $\C$ is a contour in $\D$ around the point $\sigma=s$,
we obtain by using the recursion \pref{i14}
\begin{equation} \label{q1}
f_n(s)=\frac{1}{2\pi i} \int_{\C} Q_1(\sigma,\la,s) f_{n-1}(\sigma)\,d\sigma,
\end{equation}
where $\C$ is a contour in $\D$ around the points $\sigma=\la$ and $\sigma=s$, 
and
\begin{equation} \label{q2}
Q_1(\sigma,\la,s)=-\frac{1}{2(\sigma-\la)}\left[Q_0
+2\sigma \frac{\partial}{\partial\sigma}Q_0\right]
=\frac{{\sigma+s}}{2(\sigma-\la)(\sigma-s)^2}.
\end{equation}
Continuing this we obtain for $n=0,1,2,\dots$
\begin{equation} \label{q3}
f_n(s)=\frac{1}{2\pi i} \int_{\C} Q_n(\sigma,\la,s) f(\sigma)\,d\sigma,
\end{equation}
where $\C$ is a contour in $\D$ around the points $\sigma=\la$ and $\sigma=s$. 
The rational functions $Q_n$ follow from the recursion relation
\begin{equation} \label{q4}
Q_n=-\frac{1}{2(\sigma-\la)}\left[Q_{n-1}
+2\sigma \frac{\partial}{\partial\sigma}Q_{n-1}\right].
\end{equation}
For example, we have
\begin{equation} \label{q5}
Q_2(\sigma,\la,s)= \frac{3\sigma^3-\la\sigma^2-\sigma s^2-\la s^2+
6\sigma^2s-6\la\sigma s}
{4(\sigma-\la)^3(\sigma-s)^3},
\end{equation}
The coefficients $f_k(\la)$ in \pref{i13} follow from \pref{q3} by 
substituting
$s=\la$.

In order to obtain a bound for $f_n(s)$ with representation in \pref{q3},
we select a special contour. We take for $\C$ the vertical line 
$\Re\sigma=-\sigma_0$, where $\sigma_0>0$. First we write the quantities 
$Q_n$ in a special form:
\begin{equation} \label{q6}
Q_1(\sigma,\la,s)=\frac{{2\la-p+2q}}{2\,p^2\,q}=
\frac{\la}{p^2\,q}-\frac{1}{2\,p\,q}+\frac{1}{p^2},
\end{equation}
where
\begin{equation} \label{q7}
p=\sigma-s,\quad q=\sigma-\la.
\end{equation}
Similarly,
\begin{equation} \label{q8}
Q_2(\sigma,\la,s)=\frac{8\la^2q+16\la q^2+8q^3+4p\la^2-4pq^2-2\la p^2-p^2q}
{4\,p^3\,q^3},
\end{equation}
which also can be written as a sum of partial fractions in which $\sigma$ 
occurs only in the denominator. 

For $\sigma\in\C$, we write $\sigma=-\sigma_0+i\tau$, $\tau\in\RR$.
For $s\ge 0$ and $\la\ge0$, we have
\begin{equation} \label{q9}
|p|\ge|\sigma|=\sqrt{\sigma_0^2+\tau^2},\quad |q|\ge |\la+\sigma_0+i\tau|\ge
\sqrt{\sigma_0^2+\tau^2}.
\end{equation} 
This gives for $Q_1$ the bound
\begin{equation} \label{q10}
|Q_1(\sigma,\la,s)|\le \frac{\la}{(\sigma_0^2+\tau^2)^{3/2}}
+\frac3{2(\sigma_0^2+\tau^2)}.
\end{equation}
We also need a bound of $f(s)$ for $s\in\C$ (we know that this function 
is bounded on $\C$, see beginning of Section \ref{sec:EE:BR}. Let 
\begin{equation} \label{q11}
M(\sigma_0,\la)=\max_{s\in\C} |f(s)|.
\end{equation}
 This gives 
for the remainder defined in \pref{i16} the upper bound
\begin{equation} \label{q12}
|R_1(a,z)|\le 
M(\sigma_0,\la)\,\frac{4\la+3\pi\sigma_0}{4\pi\sigma_0^2}.
\end{equation}
Similarly,
\begin{equation} \label{q13}
|Q_2(\sigma,\la,s)|\le
\frac{3\la^2}{2(\sigma_0^2+\tau^2)^{5/2}}+
\frac{11\la}{4(\sigma_0^2+\tau^2)^{2}}+
\frac{13}{8(\sigma_0^2+\tau^2)^{3/2}},
\end{equation}
\begin{equation} \label{q14}
|R_2(a,z)|\le 
M(\sigma_0,\la)\,\frac{9\la^2+11\pi\la\sigma_0+26\sigma_0^2}{16\pi\sigma_0^4},
\end{equation}

\begin{equation} \label{q15}
|Q_3(\sigma,\la,s)|\le
\frac{15\la^3}{(\sigma_0^2+\tau^2)^{7/2}}+
\frac{61\la^2}{2(\sigma_0^2+\tau^2)^{3}}+
\frac{43\la}{2(\sigma_0^2+\tau^2)^{5/2}}+
\frac{81}{8(\sigma_0^2+\tau^2)^{2}},
\end{equation}
\begin{equation} \label{q16}
|R_3(a,z)|\le 
M(\sigma_0,\la)\,\frac{768\la^3+549\pi\la^2\sigma_0+
1376\la\sigma_0^2+243\pi\sigma_0^3}
{ 96\pi\sigma_0^6}.
\end{equation}

The singularity $s_-$ estimated in \pref{sm} is of order $\bo(\la)$, and 
it follows that $\sigma_0$ can also be taken of order  $\bo(\la)$. When we 
choose $\sigma_0$ in this way,  
the quantities $Q_n$, for this choice of $\C$,  
are of order $\bo(1/\la^{n+1})$ as $\la\to\infty$, and the remainders are 
of order $\bo(1/\la^{n})$. The coefficients $f_k(\la)$ also are of
order $\bo(1/\la^{n})$. This follows from \pref{i144} and properties of 
$\phi_k(\wt \tau)$.

In the above method the estimates of higher $Q_n$ are quite easy to obtain;
only the estimate of $|f(s)|$ is needed for obtaining 
bounds of $R_n(z,a)$. From numerical verifications 
we conclude that $|f(s)|$ is maximal on $\C$ at $s=-\sigma_0$,
that is, we can take $M(\sigma_0,\la)=|f(-\sigma_0)|$.
Also we conclude that the bounds for the remainders obtained with this method 
are less realistic than those obtained by the other methods, unless $\la $ 
is becomes large.

\section{Numerical aspects}
For computing the numerical upper bounds of the
remainders in the expansion we used computer algebra for 
manipulating the formulas. This became already quite complicated,
although we didn't consider complex parameters so far. It was necessary to
develop new  algorithms for the parabolic cylinder functions, because
for certain cases we needed accurate values of the asymptotic expansions 
and the function $U(a,x)$ for computing values of the remainders in the 
expansion. 

For example, for computing $\wt R_{n}(\mu,t)$ of \pref{eq:UE.D8} for $n=3$, 
$a=100$ and $t=50$ (which corresponds with $x=1000$), the values of 
$\wt F_{\mu}(t)$ and the asymptotic series in  \pref{eq:UE.D8} (with $n=3$) 
are
$$0.99999962523819834461, \quad 0.99999962523819834799,$$
respectively, with 17 corresponding digits.
We have computed these values with Digits = 30  (in Maple), and 
developed algorithms with
adjustable precision, based on quadrature methods. We will present these 
Maple codes for computing the function $U(a,x)$ for real parameters
on the web site of our project 
({\tt http://turing.wins.uva.nl/~thk/specfun/compalg.html}). 
Fortran versions of these codes 
will be developed in a different project.
\newpage
\bibliographystyle{plain}
\bibliography{pcf}
\end{document}